\documentclass[11pt]{article}
\usepackage{amssymb}
\newcommand{\qed}{\rule{1ex}{1ex}}
\newcommand{\proof}{{\sc Proof}}
\newtheorem{definition}{Definition}[section]
\newtheorem{lemma}[definition]{Lemma}
\newtheorem{theorem}[definition]{Theorem}
\newtheorem{conjecture}[definition]{Conjecture}
\newtheorem{proposition}[definition]{Proposition}
\newtheorem{corollary}[definition]{Corollary}
\newtheorem{example}[definition]{{\sc Example}}

\newtheorem{claim}[definition]{Claim}
\newtheorem{fact}[definition]{Fact}

\begin{document}

\title{Finite models, stability, and Ramsey's theorem}
\author{Doug Ensley \\
Department of Math \\
and Computer Science \\
Shippensburg University \\
Shippensburg, PA 17257 \and Rami Grossberg \\
Department of \\
Mathematical Sciences \\
Carnegie Mellon University \\
Pittsburgh, PA 15213}
\date{\today}
\maketitle

\begin{abstract}
We prove some results on the border of Ramsey theory (finite partition
calculus) and model theory. Also a beginning of classification theory for
classes of finite models is undertaken.
\end{abstract}

\section{Introduction}

Frank Ramsey in his fundamental paper (see \cite{ra} and pages 18-27 of \cite
{grs}) was interested in ``a problem of formal logic.'' He proved the result
now known as ``(finite) Ramsey's theorem'' which essentially states

\begin{quotation}
For all $k, r, c < \omega$, there is an $n < \omega$ such that however the $r
$ -- subsets of $\{ 1, 2, \ldots, n \}$ are $c$ -- colored, there will exist
a $k$ -- element subset of $\{ 1, 2, \ldots, n \}$ which has all its $r$ --
subsets the same color.
\end{quotation}

(We will let $n(k,r,c)$ denote the smallest such $n$.) Ramsey proved this
theorem in order to construct a finite model for a given finite universal
theory so that the universe of the model is canonical with respect to the
relations in the language. (For model theorists ``canonical'' means $\Delta $
-- indiscernible as in Definition \ref{def1}).

Much is known about the order of magnitude of the function $n(k,r,c)$ and
some of its generalizations (see \cite{ehmr}, for example). An upper bound
on $n(k, r, c)$ is an $(r-1)$ -- times iterated exponential of a polynomial
in $k$ and $c$. Many feel that the upper bound is tight. However especially
for $r \geq 3$ the gap between the best known lower and upper bounds is huge.

In 1956 A. Ehrenfeucht and A. Mostowski \cite{ehmo} rediscovered the
usefulness of Ramsey's theorem in logic and introduced the notion we now
call indiscernibility. Several people continued exploiting the connections
between partition theorems and logic (i.e. model theory), among them M.
Morley (see \cite{mo1} and \cite{mo2}) and S. Shelah who has published a
virtually uncountable number of papers related to indiscernibles (see \cite
{sh}). Morley \cite{mo2} used indiscernibles to construct models of very
large cardinality (relative to the cardinality of the reals) ---
specifically, he proved that the Hanf number of $L_{\omega _1,\omega }$ is 
$\beth _{\omega _1}$.

One of the most important developments in mathematical logic --- certainly
the most important in model theory --- in the last 30 years is what is known
as ``classification theory'' or ``stability theory''. There are several
books dedicated entirely to some aspects of the subject, including books by
J. Baldwin \cite{bal}, D. Lascar \cite{las}, S. Shelah \cite{sh}, and A.
Pillay \cite{pi}. Lately Shelah and others have done extensive work in
extending classification theory from the context of first order logic, to
the classification of arbitrary classes of models, usually for infinitary
logics extending first order logic (for example see \cite{BlSh:330}, \cite
{BlSh:360}, \cite{BlSh:393}, \cite{gr2}, \cite{[GrSh 266]}, \cite{MaSh}, 
\cite{sh4}, \cite{shh}, \cite{Sh:299}). \cite{Sh:299} contains several
philosophical and personal comments about why this research is interesting,
and \cite{shtape} is the video tape of Shelah's plenary talk at the
International Congress of Mathematics at Berkeley in 1986.

This raises a question of fundamental importance: \textit{Is there a
classification theory for finite structures?} In a more philosophical
context: \textit{Is the beautiful classification theory of Shelah completely
detached >from finite mathematics?} One of the fundamental difficulties to
developing model theory for finite structures is the choice of an
appropriate ``submodel'' relation --- in category-theoretic terminology, the
choice of a natural morphism. In classification theory for elementary
classes (models of a first order, usually complete theory) the right notion
of morphism is ``elementary embedding'' defined using the relation $M \prec N
$, which is a strengthening of the notion of submodel (denoted by $M
\subseteq N$). Unfortunately for finite structures, $M \prec N$ always
implies $M=N$. Moreover, in many cases even $M \subseteq N$ implies $M=N$
(e.g., when $N$ is a group of prime order). We need a substitute.

One of the basic observations to make is that when we limit our attention to
structures in a relational language only (i.e., no function symbols), then 
$M\subseteq N$ does not imply $M=N$. In general this seems to be insufficient
to force the substructure to inherit some of the properties of the bigger
structure. It was observed already by Ramsey (in \cite{ra}) that if 
$M\subseteq N$, then for every universal sentence $\phi $, $N\models \phi $
implies $M \models \phi $. So when studying the class of models of a
universal first order theory, the relation $M \subseteq N$ is reasonable,
but it is not for more complicated theories (e.g. not every subfield of an
algebraically closed field is algebraically closed). Such a concept for
classes of finite structures is introduced in Definition \ref{submodel}.

This paper has several goals:

\begin{enumerate}
\item  \textit{Study Ramsey numbers for definable coloring inside models of
a stable theory. }

This can be viewed as a direct extension of Ramsey's work, namely by taking
into account the first order properties of the structures. A typical example
is the field of complex numbers $\langle \mathbf{C,+,\cdot \rangle }$. It is
well known that its first order theory $Th(\mathbf{C)}$ has many nice
properties --- it is $\aleph _1$ -- categorical and thus is $\aleph _0$ --
stable and has neither the order property nor the finite cover property. We
will be interested in the following general situation.

Given a first order (complete) theory $T$, and (an infinite) model $M\models
T$. Let $k,r,$ and $c$ be natural numbers, and let $F$ be a coloring of a
set of $r$ -- tuples from $M$ by $c$ colors which is definable by a first
order formula in the language $L(T)$ (maybe with parameters from $M$). Let $
n \stackrel{\mathrm{def}}{=} n_F(k,r,c)$ be the least natural number such
that for every $S \subseteq |M|$ of cardinality $n$, if $F:[S]^r \rightarrow
c$ then there exists $S^{*}\subseteq S$ of cardinality $k$ such that $F$ is
constant on $[S^{*}]^r$. It turns out that for stable theories, (or even for
theories without the independence property) we get better upper bounds than
for the general Ramsey numbers. This indicates that one can not improve the
lower bounds by looking at stable structures.

\item  \textit{Introduce stability-like properties ( e.g. }$\mathit{n}$
\textit{-order property, }$\mathit{k}$\textit{-independence property, }$
\mathit{d}$\textit{-cover property), as well as averages of finite sequences
of indiscernibles.}

Some of the interconnections and the effect on the existence of
indiscernibles are presented.

\item  \textit{Develop classification theory for classes of finite
structures. In particular introduce a notion that correspond to stable
amalgamation, and show that it is symmetric for many models. }

See Example \ref{ex.fields}.

\item  \textit{Bring down uncountable techniques to a finite context. }

We believe that much of the machinery developed (mainly by Shelah) to deal
with problems concerning categoricity of infinitary logics and the behavior
of the spectrum function at cardinalities $\geq \beth _{\omega _1}$ depends
on some very powerful combinatorial ideas. We try here to extract some of
these ideas and present them in a finite context.
\end{enumerate}

Shelah \cite{sh} proved that instability is equivalent to the presence of
either the strict order property or the independence property. In a
combinatorial setting, stability implies that for arbitrarily large sets,
the number of types over a set is polynomial in the cardinality of the set.
We address the finite case here in which we restrict our attention to when
the number of $\phi $- types over a finite set is bounded by a polynomial in
the size of the set of parameters.

First we find precisely the degree of the polynomial bound on the number of
these types given to us by the absence of the strict order or independence
properties. This is an example of something relevant in the finite case
which is of no concern in the usual classification theory framework.

Once we have these sharper bounds we can find sequences of indiscernibles in
the spirit of \cite{sh}. It should be noted here that everything we do is
``local'', involving just a single formula (or equivalently a finite set of
formulas). We then work through the calculations for uniform hypergraphs as
a case study. This raises questions about ``stable'' graphs and hypergraphs
which we begin to answer.

In the second half of the paper, we examine classes of finite structures in
the framework of Shelah's classification for non-elementary classes (see 
\cite{sh3}). In particular, we make an analogy to Shelah's ``abstract
elementary classes'' and prove results similar to his.

\noindent \textbf{Notation}: Everything is standard. We will typically treat
natural numbers as ordinals (i.e., $n = \{0, 1, \ldots, n-1 \}$). Often $x$, 
$y$, and $z$ will denote free variables, or finite sequences of variables
--- it should be clear from the context whether we are dealing with
variables or with sequences of variables. When $x$ is a sequence, we let $
l(x)$ denote its length. $L$ will denote a similarity type (a.k.a. language
or signature), $\Delta $ will stand for a finite set of $L$ formulas. $M$
and $N$ will stand for $L$ - structures, $|M|$ the universe of the structure 
$M$, and $\Vert M\Vert $ the cardinality of the universe of $M$. Given a
fixed structure $M$, subsets of its universe will be denoted by $A$, $B$, $C$
, and $D$. So when we write $A \subseteq M$ we really mean that $A\subseteq
|M|$, while $N\subseteq M$ stands for ``$N$ is a submodel of $M$''. When $M$
is a structure then by $a \in M$ we mean $a \in |M|$, and when $a$ is a
finite sequence of elements, then $a\in M$ stands for ``all the elements of
the sequence $a$ are elements of $|M|$''.

Since all of our work will be inside a given structure $M$ (with the
exception of Section \ref{graphsec}), all the notions are relative to it.
For example for $a\in M$ and $A\subseteq M$ we denote by $tp_\Delta (a,A)$
the type $tp_\Delta (a,A,M)$ which is $\{\phi (x;b):M\models \phi [a;b],b\in
A,\phi (x;y)\in \Delta \}$ and if $A\subseteq M$ then $S_\Delta (A,M)
\stackrel{\mathrm{def}}{=}\{{tp_\Delta (a,A):a\in M\}}$. Note that in \cite
{sh} $S_\Delta (A,M)$ denotes the set of all complete $\Delta $ -types with
parameters from $A$ that are consistent with $Th(\langle M,c_a\rangle _{a\in
A}$. It is important for us to limit attention to the types realized in $M$
in order to avoid dependence on the compactness theorem. It is usually
important that $\Delta $ is closed under negation, so when $\Delta =\{\phi
,\neg \phi \}$, instead of writing $tp_\Delta (\cdots )$ and $S_\Delta
(\cdots )$ we will write $tp_\phi (\cdots )$ and $S_\phi (\cdots )$,
respectively.

\section{The effect of the order and independence properties on the number
of local types}

In this section, we fix some notation and terms and then define the first
important concepts. In the following definition, the first three parts are
>from \cite{sh}, $(4)$ is a generalization of a definition of Shelah, and $
(5)$ is from Grossberg and Shelah \cite{GrSh}.

\begin{definition}
\label{def1}

\begin{enumerate}
\item  For a set $\Delta $ of $L$ -- formulas and a natural number $n$, a 
\underline{$(\Delta ,n)$ -- type over a set} $A$ is a set of formulas of the
form $\phi ({x};a)$ where $\phi (x;y)\in \Delta $ and $a\in A$ with $l(x)=n$
. If $\Delta =L$, we omit it, and we just say ``$\phi $ -- type'' for a $
(\{\phi (x;y),\neg \phi (x;y)\},l(x))$ -- type.

\item  Given a $(\Delta ,n)$ -- type $p$ over $A$, define $dom(p)=\{a\in
A\;:\;$ for some $\phi \in \Delta ,\,\phi (x;a)\in p\}$.

\item  A type $p$ \underline{$(\Delta _0,\Delta _1)$ -- splits over} $
B\subseteq dom(p)$ if there is a $\phi (x;y)\in \Delta _0$ and $b,c\in dom(p)
$ such that $tp_{\Delta _1}(b,B)=tp_{\Delta _1}(c,B)$ and $\phi (x;b),\neg
\phi (x;c)\in p$. If $p$ is a $\Delta $ -- type and $\Delta _0=\Delta
_1=\Delta $, then we just say $p$ \underline{splits over $B$}.

\item  \label{def1d} We say that $(M,\phi (x;y))$ has \underline{the $k$ -
independence property} if there are $\{a_i:i<k\}\subseteq M$, and $
\{b_w:w\subseteq k\}\subseteq M$, such that $M\models \phi [a_i;b_w]$ if and
only if $i\in w$. We will say that \underline{$M$ has the $k$ --
independence property} when there is a formula $\phi $ such that $(M,\phi )$
does.

\item  $(M,\phi (x;y))$ has \underline{the $n$ -- order property} (where $
l(x)=l(y)=k$) if there exists a set of $k$ -- tuples $\{a_i:i<n\}\subseteq M$
such that $i<j$ if and only if $M\models \phi [a_i,a_j]$ for all $i,j<n$. We
will say that $M$ has \underline{the $n$ -- order property} if there is a
formula $\phi$ so that $(M,\phi )$ has the $n$ -- order property.
\end{enumerate}
\end{definition}

\noindent \textsc{Warning:} This use of ``order property'' corresponds to
neither the order property nor the strict order property in \cite{sh}. The
definition comes rather from \cite{gr}.

The following monotonicity property is immediate from the definitions.

\begin{proposition}
Let sets $B\subseteq C\subseteq A$ and a complete $(\Delta, n)$ -- type $p$
be given with $Dom(p)\subseteq A$. If $p$ does not split over $B$, then $p$
does not split over $C$.
\end{proposition}

\begin{fact}
\label{unstable}( Shelah see \cite{sh}) Let $T$ be a complete first order
theory. The following conditions are equivalent:

\begin{enumerate}
\item  $T$ is unstable.

\item  There are $\phi (x;y)\in L(M)$, $M\models T$, and $\{a_n:n<\omega
\}\subseteq M$ such that $l(x)=l(y)=l(a_n)$, and for every $n,k<\omega $ we
have \\$n<k\Leftrightarrow M\models \phi [a_n;a_k]$.
\end{enumerate}
\end{fact}

The compactness theorem gives us the following

\begin{corollary}
\label{finitecor} Let $T$ be a stable theory, and suppose that $M\models T$
is an infinite model.

\begin{enumerate}
\item  For every $\phi (x,y)\in L(M)$ there exists a natural number $n_{\phi}
$ such that $(M,\phi )$ does not have the $n_{\phi}$ -- order property.

\item  For every $\phi (x,y)\in L(M)$ there exists a natural number $k_{\phi}
$ such that $(M,\phi )$ does not have the $k_{\phi}$ -- independence
property.

\item  If $T$ is categorical in some cardinality greater than $|T|$, then
for every $\phi (x,y)\in L(M)$ there exists a natural number $d_{\phi}$ such
that $(M,\phi )$ does not have the $d_{\phi}$ -- cover property (see
Definition \ref{coverp}).
\end{enumerate}
\end{corollary}

We first establish that the failure of either the independence property or
the order property for $\phi$ implies that there is a polynomial bound on
the number of $\phi$ -- types. The more complicated of these to deal with is
the failure of the order property. At the same time this is perhaps the more
natural property to look for in a given structure. The bounds in this case
are given in Theorem \ref{thm2}. The failure of the independence property
gives us a far better bound (i.e., smaller degree polynomial) with less
work. Theorem \ref{nip_poly} reproduces this result of Shelah paying
attention to the specific connection between the bound and where the
independence property fails.

This first lemma is a finite version of Lemma 5 from \cite{gr}.

\begin{lemma}
\label{lem1} Let $\phi (x;y)$ be a formula in $L$, $n$ a positive integer, $
s=l(y)$, $r=l(x)$, $\psi (y;x)=\phi (x;y)$. Suppose that $\{A_i\subseteq M
\;:\;i\leq 2n\}$ is an increasing chain of sets such that for every $
B\subseteq A_i$ with $|B|\leq 3sn$, every type in $S_\phi (B,M)$ is realized
in $A_{i+1} $. Then if there is a type $p\in S_\phi (A_{2n},M)$ such that
for all $i < 2n$, $p|A_{i+1}$ $(\psi ,\phi )$ -- splits over every subset of 
$A_i$ of size at most $3sn$, then $(M,\rho )$ has the $n-$ order property,
where 
\[
\rho (x_0,x_1,x_2;y_0,y_1,y_2)\stackrel{\mathrm{def}}{=}[\phi
(x_0;y_1)\leftrightarrow \phi (x_0;y_2)] 
\]
\end{lemma}

{\proof} Choose $d \in M$ realizing $p$. Define $\{a_i,b_i,c_i\in
A_{2i+2}\;:\;i<n\}$ by induction on $i$. Assume for $j<n$ that we have
defined these for all $i<j$. Let $B_j=\bigcup \{a_i,b_i,c_i\;:\;i<j\}$.
Notice that $|B_j|\leq 3sj<3sn$, so by the assumption, $p|A_{2j+1}$ $(\psi
,\phi )$ -- splits over $B_j$. That is, there are $a_j$, $b_j\in A_{2j+1}$
such that 
\[
tp_\psi (a_j,B_j,M)=tp_\psi (b_j,B_j,M), 
\]
and 
\[
M\models \phi [d;a_j]\wedge \neg \phi [d;b_j]. 
\]
Now choose $c_j\in A_{2j+1}$ realizing $tp(d,B_j\cup \{ a_j,b_j\},M)$ (which
can be done since $|B_j\cup\{a_j,b_j\}|\leq 3sj+2s<3s(j+1)\leq 3sn$). This
completes the inductive definition.

For each $i$, let $d_i=c_ia_ib_i$. We will check that the sequence of $d_i$
and the formula 
\[
\rho (x_0,x_1,x_2;y_0,y_1,y_2)\stackrel{\mathrm{def}}{=}[\phi
(x_0;y_1)\leftrightarrow \phi (x_0;y_2)] 
\]
witness the $n$ -- order property for $M$.

If $i<j<n$, then $c_i\in B_j$. By choice of $a_j$ and $b_j$, $tp_\psi
(a_j,B_j,M)=tp_\psi (b_j,B_j,M)$, so in particular, 
\[
M\models \phi [c_i;a_j]\leftrightarrow \phi [c_i;b_j] 
\]
That is, $M\models \rho [d_i;d_j]$.

On the other hand, if $i\leq j<n$, then $\phi (x;a_i)\in tp_\phi (d,B_j \cup
\{a_j,b_j\},M)$ and $\phi (x;b_i) \not{\in }tp_\phi (d,B_j \cup \{a_j,
b_j\},M)$, and so, by the  choice of $c_j$, we have that 
\[
M\models \phi [c_j;a_i]\wedge \neg \phi [c_j;b_i]. 
\]
That is, $M\models \neg \rho [d_j;d_i]$ in this case. {\qed}

In order to see the relationship between this definition of the order
property and Shelah's, we mention Corollary \ref{cor1} below. Note that it
is the formula $\phi$, not the $\rho$ of Lemma \ref{lem1}, which has the
weak order property in the Corollary.

\begin{definition}
$(M,\phi )$ has \underline{the weak $m$ -- order property} if there exist $
\{d_i:i<m\}\subseteq M$ such that for each $j<m$, 
\[
M\models \exists x\bigwedge_{i<m}\phi (x;d_i)^{\mbox{if}(i\geq j)} 
\]
\end{definition}

\noindent\textsc{Remark:} This is what Shelah \cite{sh} calls the $m$ --
order property.

\begin{definition}
\label{arrowdef} We write $x\rightarrow (y)_b^a$ if for every partition $\Pi 
$ of the $a$ - element subsets of $\{1,\ldots ,x\}$ with $b$ parts, there is
a $y$ - element subset of $\{1,\ldots ,x\}$ with all of \emph{its} $a$ --
element subsets in the same part of $\Pi $.
\end{definition}

\begin{corollary}
\label{cor1}

\begin{enumerate}
\item  If in addition to the hypotheses of Lemma \ref{lem1} we have that $
(2n)\rightarrow (m+1)_2^2$, then $\phi$ has the weak $m$ -- order property
in $M$.

\item  If in addition to the hypotheses of Lemma \ref{lem1} we have that $
n\geq \frac{2^{2m-1}}{\pi m}$, then $\phi $ has the weak $m$ -- order
property in $M$.
\end{enumerate}
\end{corollary}

{\proof} (This is essentially \cite{sh} I.2.10(2))

\begin{enumerate}
\item  Let $a_i$, $b_i$, $c_i$ for $i<n$ be as in the proof of Lemma \ref
{lem1}. For each pair $i<j\leq n$, define 
\[
\chi ({i,j}):=\left\{ 
\begin{array}{ll}
1 & \mbox{ if $M \models \phi[c_i;  a_j]$} \\ 
0 & \mbox{otherwise.}
\end{array}
\right. 
\]

Since $(2n)\rightarrow (m+1)_2^2$, we can find a subset $I$ of $2n$ of
cardinality $m+1$ on which $\chi $ is constant and which we can enumerate as 
$I = \{i_0,\ldots ,i_m\}$.

If $\chi $ is 1 on $I$, then for every $k$ with $1\leq k\leq m+1$ 
\[
\{\neg \phi (x;b_{i_j})^{\mbox{\em if}(j>k)}\;:\;1\leq j<m\} 
\]
is realized by $c_{i_{k-1}}$. Therefore, the sequence $\{b_{i_0},\ldots
,b_{i_m}\}$ witnesses the weak $m$ -- order property of $\phi $ in $M$.

On the other hand, if $\chi $ is 0 on $I$, then for every $k$ with $1\leq
k\leq m+1$ 
\[
\{\neg \phi (x;a_{i_j})^{\mbox{\em if}(j>k)}\;:\;1\leq j<m\} 
\]
is realized by $c_{i_{k-1}}$. Therefore, the sequence $\{a_{i_0},\ldots
,a_{i_m}\}$ witnesses the weak $m$ -- order property of $\neg \phi $ in $M$.
Of course, it is equivalent for $\phi $ and $\neg \phi $ to have the weak $m$
-- order property in $M$.

\item  By Stirling's formula, $n\geq \frac{2^{2m-1}}{\pi m}$ implies that $
n\geq \frac 12{\left( 
\begin{array}{c}
{2m} \\ 
{m}
\end{array}
\right) }$, and from \cite{grs}, $n\geq \frac 12{\left( 
\begin{array}{c}
{2m} \\ 
{m}
\end{array}
\right) }$ implies that $(2n)\rightarrow (m+1)_2^2$.
\end{enumerate}

{\qed}

We can now establish the relationship between the number of types and the
order property.

\begin{theorem}
\label{thm2} If $\phi (x;y)\in L(M)$ is such that 
\[
\rho (x_0,x_1,x_2;y_0,y_1,y_2)\stackrel{\mathrm{def}}{=}\phi
[(x_0;y_1)\leftrightarrow \phi (x_0;y_2)] 
\]
does not have the $n$ -- order property in $M$, then for every set $
A\subseteq M$ with $|A|\geq 2$, we have that $|S_\phi (A,M)|\leq 2n|A|^k$,
where $k=2^{(3ns)^{t+1}}$ for $r=l(x)$, and $s=l(y)$, and $t=\max \{r,s\}$.
\end{theorem}

{\proof} Suppose that there is some $A\subseteq M$ with $|A|\geq 2$ so that $
|S_\phi (A,M)|>(2n)|A|^k$. Let $\psi (y;x)=\phi (x;y)$, $m=|A|$, and let $
\{a_i\;:\;i\leq (2n)m^k\}\subseteq M$ be witnesses to the fact that $|S_\phi
(A,M)|>(2n)m^k$. (That is, each of these tuples realizes a different $\phi $
-- type over $A$.) Define $\{A_i:i<2n\}$, satisfying

\begin{enumerate}
\item  $A\subseteq A_i\subseteq A_{i+1} \subseteq M$,

\item  $|A_i|\leq c^{e(i)}m^{(3ns)^i}$,where $c:=2^{2+(3sn)^t}$ and $e(i):= 
\frac{(3ns)^i-1}{3ns-1}$, and

\item  for every $B\subseteq A_i$ with $|B|\leq 3sn$, every $p\in S_\phi
(B,M)\bigcup {S_\psi (B,M)}$ is realized in $A_{i+1}$.
\end{enumerate}

To see that this can be done, we need only check the cardinality
constraints. There are at most $|A_i|^{3sn}$ subsets of $A_i$ with
cardinality at most $3sn$, and over each such subset $B$, there are at most $
2^{(3sn)^r}$ and $2^{(3sn)^s}$ types in $S_\psi (B,M)$ and $S_\phi (B,M)$,
respectively, so there are at most $2^{(3sn)^r}+2^{(3sn)^s}\leq
2^{1+(3sn)^t} $ types in $S_\psi (B,M)\bigcup {S_\phi (B,M)}$ for each such $
B$. Therefore, $A_{i+1}$ can be defined so that 
\begin{eqnarray*}
|A_{i+1}| & \leq & |A_i|+(2^{1+(3sn)^t})|A_i|^{3sn} \\
& \leq & c|A_i|^{3sn} \\
& \leq & c(c^{e(i)}m^{(3ns)^i})^{3sn} \\
& = & c^{1+e(i)(3sn)}m^{(3sn)^{i+1}} \\
& = & c^{e(i+1)}m^{(3sn)^{i+1}}
\end{eqnarray*}

\begin{claim}
\label{thm2c1} There is a $j<(2n)m^k$ such that for every $i<2n$ and every $
B\subseteq A_i$ with $|B|\leq 3sn$, $tp(a_j,A_{i+1})$ $(\psi ,\phi )$ --
splits over $B$.
\end{claim}

{\proof} (Of Claim \ref{thm2c1}) Suppose not. That is, for every $j\leq
(2n)m^k$, there is an $i(j)<2n$ and a $B\subseteq A_{i(j)}$ with $|B|\leq 3sn
$, so that $tp(a_j,A_{i(j)+1})$ does not $(\psi ,\phi )$ -- split over $B$.
Since $i$ is a function from $1+(2n)m^k$ to $2n$, there must be a subset $S$
of $1+(2n)m^k$ with $|S|>m^k$, and an integer $i_0<2n$ such that for all $
j\in S$, $i(j)=i_0$. Now similarly, there are less than $|A_{i_0}|^{3sn}$
subsets of $A_{i_0}$, with cardinality at most $3sn$, so there is a $
T\subseteq S$ with 
\[
|T|>\frac{m^k}{|A_{i_0}|^{3sn}} 
\]
and a $B_0\subseteq A_{i_0}$, with $|B_0|\leq 3sn$ such that for all $j\in T$
, $tp(a_j,A_{i_0+1})$ does not $(\psi ,\phi )$ -- split over $B_0$. Since $
|A_{i_0}|\leq c^{e(i_0)}m^{(3ns)^{i_0}}\leq (cm)^{(3sn)^{2n}}$, then 
\[
|T|\geq \frac{m^k}{(cm)^{(3sn)^{2n}}} 
\]
Let $C\subseteq A_{i_0+1}$ be obtained by adding to $B_0$, realizations of
every type in $S_\phi (B_0,M)\bigcup {\ S_\psi (B_0,M)}$. This can clearly
be done so that $|C|\leq 3ns+2^{(3ns)^r}+2^{(3ns)^s}$. The maximum number of 
$\phi $ -- types over $C$ is at most $2^{|C|^s}\leq 2^{c^s}$.

\begin{claim}
\label{claim11} \label{thm2c2} $m^{k-(3ns)^{2n}}>(2^{c^s})(c^{(3ns)^{2n}})$
\end{claim}

{\proof} (Of Claim \ref{thm2c2}) Since $c=2^{2+(3ns)^t}$, we have $
c^s+(3ns)^{2n}(2+(3ns)^t)$ as the exponent on the right-hand side above.
Since $m\geq 2$, it is enough to show that 
\begin{eqnarray*}
k>(c^s+(3ns)^{2n}(2+(3ns)^t)+(3ns)^{2n} \\
=2^{s(2+(3ns)^t)}+(3ns)^{2n}(3+(3ns)^t)
\end{eqnarray*}
This follows from the definition of $k$ (recall that $k=2^{(3ns)^{t+1}}$),
so we have established Claim \ref{thm2c2}. {\qed}$_{\ref{claim11}}$

Therefore, $|T|$ is greater than the number of $\phi $ -- types over $C$, so
there must be $i\neq j\in T$ such that $tp_\phi (a_i,C)=tp_\phi (a_j,C)$.
Since $tp_\phi (a_i,A)\neq tp_\phi (a_j,A)$, we may choose $a\in A$ so that $
M\models \phi [a_i,a]\wedge \neg \phi [a_j,a]$. Now choose $a^{\prime }\in C$
so that $tp_\psi (a,B_0)=tp_\psi (a^{\prime },B_0)$ (this is how $C$ is
defined after all). Since $tp_\phi (a_i,A_{i_0+1})$ does not $(\psi ,\phi )$
-- split over $B_0$, we have that 
\[
\phi (x;a)\in tp_\phi (a_i,A_{i_0+1})\mbox{ if and only if }\phi
(x;a^{\prime })\in tp_\phi (a_i,A_{i_0+1}) 
\]
so $M\models \phi [a_i,a^{\prime }]\wedge \neg \phi [a_j,a^{\prime }]$,
contradicting the fact that $tp_\phi (a_i,C)=tp_\phi (a_j,C)$ and thus
completing the proof of Claim \ref{thm2c1}. Now letting $j$ be as in Claim 
\ref{thm2c1} above and applying Lemma \ref{lem1} completes the proof of
Theorem \ref{thm2}. {\qed}

Theorem \ref{thm4} below gives a better result under different assumptions.
The next lemma is II, 4.10, (4) in \cite{sh}. It comes from a question due
to Erd\H os about the so-called ``trace'' of a set system which was answered
by Shelah and Perles \cite{sh2} in 1972. Purely combinatorial proofs (i.e.,
proofs in the language of combinatorics) can also be found in most books on
extremal set systems (e.g., Bollobas \cite{bo}).

\begin{lemma}
\label{lem3} If $S$ is any family of subsets of the finite set $I$ with 
\[
|S|>\sum_{i<k}{\left( 
\begin{array}{c}
{|I|} \\ 
{i}
\end{array}
\right) } 
\]
then there exist $\alpha _i\in I$ for $i<k$ such that for every $w\subseteq
k $ there is an $A_w\in S$ so that $i\in w\Leftrightarrow \alpha _i\in A_w$.
(The conclusion here is equivalent to $trace(I)\geq k$ in the language of 
\cite{bo}.)
\end{lemma}

{\proof} See Theorem 1 in Section 17 of \cite{bo} or Ap.1.7(2) in \cite{sh}. 
{\qed}

\begin{theorem}
\label{thm4}\label{nip_poly} If $\phi (x;y)\in L(M)$ ($r=l(x)$, $s=l(y)$)
does not have the $k$ -- independence property in $M$, then for every set $
A\subseteq M$, if $|A|\geq 2$, then $|S_\phi (A,M)|\leq |A|^{s(k-1)}$.
\end{theorem}

{\proof} (Essentially \cite{sh}, II.4.10(4)) Let $F$ be the set of $\phi $
-- formulas over $A$. Then 
\[
|F|<|A|^s. 
\]
So if $|S_\phi (A,M)|>|A|^{s(k-1)}$, then certainly 
\[
|S_\phi (A,M)|>\sum_{i<k}\left( 
\begin{array}{c}
{|F|} \\ 
{i}
\end{array}
\right) , 
\]
in which case Lemma \ref{lem3} can be applied to $F$ and $S_\phi (A,M)$ to
get witnesses to the $k$ -- independence property in $M$, a contradiction. {
\qed}

The ``moral'' of Theorem \ref{thm2} and Theorem \ref{thm4} is that when $\phi
$ has some nice properties, there is a bound on the number of $\phi$ --
types over $A$ which is polynomial in $\vert A \vert$. Note that the
difference between the two properties is that the degree of the polynomial
in the absence of the $k$ -- independence property is linear in $k$ while in
the absence of the $n$ -- order property the degree is exponential in $n$.
Also the bounds on $\phi$ -- types in the latter case hold when a formula $
\rho$ related to $\phi$ (as opposed to $\phi$ itself) is without the $n$ --
order property.

Another property discovered by Keisler (in order to study saturation of
ultrapowers, see \cite{ke}), and studied extensively by Shelah is the
``finite cover property'' (see \cite{sh}) whose failure essentially provides
us with a strengthening of the compactness theorem.

\begin{definition}
\label{coverp} We say that $(M,\phi )$ does not have \underline{the $d$ -
cover property} if for every $n\geq d$ and $\{b_i:i<n\}\subseteq M$, \textbf{
if} 
\[
\left( \forall w\subseteq n \left[ |w|<d\Rightarrow M\models \exists
x\bigwedge_{i\in w}\phi (x;b_i)\right] \right) 
\]
\textbf{then} 
\[
M\models \exists x \bigwedge_{i<n}\phi (x;b_i). 
\]
\end{definition}

\begin{example}
If $M=(M,R)$ is the countable random graph, then $(M,R)$ fails to have the 2
-- cover property. If $M$ is the countable universal homogeneous
triangle-free graph, then $(M,R)$ fails to have the 3 -- cover property.
\end{example}

\section{Indiscernible sequences in large finite sets}

\noindent\textsc{Note:} The next definition is an interpolant of Shelah's 
\cite{sh}, I.2.3, and Ramsey's notion of canonical sequence.

\begin{definition}
\begin{enumerate}
\item  A sequence $I=\langle a_i\;:\;i<n\rangle \subseteq M$ is called a $
(\Delta ,m)$ -- \underline{indiscernible sequence over} $A\subseteq M$
(where $\Delta $ is a set of $L(M)$ -- formulas) if for every $i_0<\ldots
<i_{m-1}\in I$, $j_0<\ldots <j_{m-1}\in I$ we have that $tp_\Delta
(a_{i_0}\cdots a_{i_{m-1}},A,M)=tp_\Delta (a_{j_0}\cdots a_{j_{m-1}},A,M)$

\item  A set $I=\{ a_i\;:\;i<n \} \subseteq M$ is called a $(\Delta ,m)$ -- 
\underline{indiscernible set over} $A\subseteq M$ if and only if for every $
\{i_0,\ldots ,i_{m-1}\}$, $\{j_0,\ldots ,j_{m-1}\}\subseteq I$ we have

$tp_\Delta (a_{i_0}\cdots a_{i_{m-1}},A,M)=tp_\Delta (a_{j_0}\cdots
a_{j_{m-1}},A,M)$.
\end{enumerate}
\end{definition}

Note that if $\phi(x; b) \in tp(a_0 \ldots a_{m-1}, B, M)$, then necessarily 
$l(x) = m \cdot l(a_0)$.

\begin{example}
\begin{enumerate}
\item  In the model $M^m_n = \langle m,0,1,\chi \rangle $ ($n\leq m<\omega $
) where $\chi $ is function from the increasing $n$ -- tuples of $m$ to $
\{0,1\}$, any increasing enumeration of a monochromatic set is an example of
a $(\Delta ,1)$ -- indiscernible sequence over $\emptyset$ where $\Delta
=\{\chi(x)=0,\chi (x)=1\}$.

\item  In a graph $(G,R)$, cliques and independent sets are examples of $
(R,2)$ -- indiscernible sets over $\emptyset $.
\end{enumerate}
\end{example}

Recall that in a stable first order theory, every sequence of indiscernibles
is a set of indiscernibles. In our finite setting this is also true if the
formula fails to have the $n$ -- order property. The argument below follows
closely that of Shelah \cite{sh}.

\begin{theorem}
\label{thm6} If $M$ does not have the $n$ -- order property, then any
sequence $I=\langle a_i\;:\; i < n + m - 1 \rangle \subseteq M$ which is $
(\phi,m)$ -- indiscernible over $B\subseteq M$ is a set of $(\phi,m)$ --
indiscernibles over $B$.
\end{theorem}

{\proof} Since any permutation of $\{1,\ldots ,n\}$ is a product of
transpositions $(k,k+1)$, and since $I$ is a $\phi$ -- indiscernible
sequence over $B$, it is enough to show that for each $b\in B$ and $k < m$, 
\[
M\models \phi [a_0\cdots a_{k-1}a_{k+1}a_k\cdots a_{m-1};b]\leftrightarrow
\phi [a_0\cdots a_{k-1}a_ka_{k+1}\cdots a_{m-1};b]. 
\]
Suppose this is not the case. Then we may choose $b\in B$ and $k<m$ so that 
\[
M\models \neg \phi [a_0\cdots a_{k-1}a_{k+1}a_k\cdots a_{m-1};b]\wedge \phi
[a_0\cdots a_{k-1}a_ka_{k+1}\cdots a_{m-1};b]. 
\]
Let $c=a_0\cdots a_{k-1}$ and $d=a_{n+k+1}\cdots a_{n+m-2}$ making $l(c)=k$
and $l(d)=m-k-2$). By the indiscernibility of $I$, 
\[
M\models \neg \phi [ca_{k+1}a_kd;b]\wedge \phi [ca_ka_{k+1}d;b]. 
\]
For each $i$ and $j$ with $k\leq i<j<n+k$, we have (again by the
indiscernibility of the sequence $I$) that 
\[
M\models \neg \phi [ca_ja_id;b]\wedge \phi [ca_ia_jd;b]. 
\]
Thus the formula $\psi (x,y;cdb)\stackrel{\mathrm{def}}{=}\phi (c,x,y,d;b)$
defines an order on $\langle a_i\;:\;k\leq i<n+k\rangle $ in $M$, a
contradiction. {\qed}

The following definition is a generalization of the notion of end-homogenous
sets in combinatorics (see section 15 of \cite{ehmr}) to the context of $
\Delta$ -- indiscernible sequences.

\begin{definition}
A sequence $I=\langle a_i\;:\;i<n\rangle \subseteq M$ is called a $(\Delta
,m)$ -- \underline{end-indiscernible sequence over} $A\subseteq M$ (where $
\Delta $ is a set of $L(M)$ -- formulas) if for every $\{i_0,\ldots
,i_{m-2}\} \subseteq n$ and $j_0,j_1<n$ both larger than $\max \{i_0,\ldots
,i_{m-2}\}$, we have 
\[
tp_\Delta (a_{i_0}\cdots a_{i_{m-2}}a_{j_0},A,M)=tp_\Delta (a_{i_0}\cdots
a_{i_{m-2}}a_{j_1},A,M) 
\]
\end{definition}

\begin{definition}
For the following lemma, let $F:\omega \rightarrow \omega $ be given, and
fix the parameters, $\alpha $, $r$, and $m$. We define the function $F^{*}$
for each $k \geq m$ as follows:

\begin{itemize}
\item  $F^{*}(0) = 1$,

\item  $F^{*}(j+1) = 1 + F^{*}(j)\cdot F(\alpha + m \cdot r \cdot j)$ for $j
< k - 2 -m$, and

\item  $F^{*}(j+1) = 1 + F^{*}(j)$ for $k-2-m \leq j<k-2$.
\end{itemize}

We will not need $j\geq k-2$.
\end{definition}

\begin{lemma}
\label{lem7} Let $\psi (x;y) = \phi (x_1,\ldots ,x_{m - 1},x_0;y)$. If for
every $B\subseteq M$, $|S_\psi (B,M)|<F(|B|)$, and $I=\{c_i\;:\;i\leq
F^{*}(k-2)\}\subseteq M$ (where $l(c_i)=l(x_i)=r$, $\alpha =|A|$), then
there is a $J\subseteq I$ such that $|J|\geq k$ and $J$ is a $(\phi,m)$ --
end-indiscernible sequence over $A$.
\end{lemma}

{\proof} (For notational convenience when we have a subset $S \subseteq I$,
we will write $\min S$ instead of the clumsier $c_{\min \{i \; : \; c_i \in
S \}}$.) We now construct $A_j=\{a_i \; : \; i\leq j\} \subseteq I$ and $S_j
\subseteq I$ by induction on $j<k-1$ so that

\begin{enumerate}
\item  \label{c71} $a_j=\min S_j$,

\item  \label{c72} $S_{j+1} \subseteq S_j$,

\item  \label{c73} $|S_j| > F^{*}(k-2-j)$, and

\item  \label{c74} whenever $\{i_0,\ldots ,i_{m-1}\}\subseteq j$ and $b\in
S_j$, 
\[
tp_\phi (a_{i_0}\cdots a_{i_{m-2}}a_j,A,M)=tp_\phi (a_{i_0}\cdots
a_{i_{m-2}}b,A,M). 
\]
\end{enumerate}

The construction is completed by taking an arbitrary $a_{k-1}\in
S_{k-2}-\{a_{k-2}\}$. (which is possible by (\ref{c73}) since $F^{*}(0)=1$),
and letting $J=\langle a_i\;:\;i<k\rangle $. We claim that $J$ will be the
desired ($\phi,m)$ -- end-indiscernible sequence over $A$.

To see this, let $\{i_0,\ldots ,i_{m-2},j_0,j_1\}\subseteq k$ with $\max
\{i_0,\ldots ,i_{m-2}\}<j_0<j_1<k$ be given. Certainly then $\{i_0,\ldots
,i_{m-2}\}\subseteq j_0$ and $a_{j_1}\in S_{j_0}$, so by (\ref{c74}) we have
that 
\[
tp_\phi (a_{i_0}\cdots a_{i_{m-2}}a_{j_0},A,M)=tp_\phi (a_{i_0}\cdots
a_{i_{m-2}}a_{j_1},A,M). 
\]

To carry out the construction, first set $a_j = c_j$ and $S_j = \{c_i \; :
\; j \leq i\leq F^{*}(k-2)\}$ for $0\leq j\leq m - 1$. Clearly, we have
satisfied all conditions in this. Now assume for some $j\geq m$ that $A_{j-1}
$ and $S_{j-1}$ have been defined satisfying the conditions.

Define the equivalence relation $\sim $ on $S_{j-1}-\{a_{j-1}\}$ by $c\sim d$
if and only if for all $\{i_0,\ldots ,i_{r-1}\}$, 
\[
tp_\phi (a_{i_1}\cdots a_{i_{m-1}}c,A,M)=tp_\phi (a_{i_1}\cdots
a_{i_{m-1}}d,A,M) 
\]
The number of $\sim $ -- classes then is at most $|S_\psi (A \cup
A_j)|<F(\alpha +m \cdot r \cdot j)$. Therefore, at least one class $S_m cdot
r cdot j$ has cardinality at least $\frac{|S_{j-1}|-1}{F(\alpha + m \cdot r
\cdot j)}$. Let $a_j=\min S_j$. By definition of $F^{*}$, $\frac{
F^{*}(k-2-j+1)}{F(\alpha + m \cdot r \cdot j)}>F^{*}(k-2-j)$, so we have
that $|S_j|>F^{*}(k-2-j)$. It is easy to see that condition (\ref{c74}) is
satisfied. {\qed}

For the following lemma, we once again need a function defined in terms of
the parameters of the problem. We will need the parameter $r$ and the
function $F^*$ defined for Lemma \ref{lem7} (which depends on $r$, $\alpha$,
and $m$). Let $f_i$ be the $F^*$ that we get when $m = i$ (and $\alpha$ and $
r$ are fixed) in Lemma \ref{lem7}.

For the following lemma define 
\[
g_i:=\left\{ 
\begin{array}{ll}
id & \mbox{if }i=0 \\ 
f_{i-1}\circ (g_{i-1}-2) & \mbox{otherwise}
\end{array}
\right. 
\]

\begin{lemma}
\label{lem8} If $J=\{a_i\;:\;i\leq g_{m-1}(k-1)\}\subseteq M$ is a $(\phi,m) 
$ -- end-indiscernible sequence over $A\subseteq M$, then there is a $
J^{\prime }\subseteq J$ such that $|J^{\prime }|\geq k$ and $J^{\prime }$ is
a $(\phi,m) $ -- indiscernible sequence over $A$.
\end{lemma}

{\proof} (By induction on $m$) Note that if $m=1$, there is nothing to do
since end - indiscernible \emph{is} indiscernible in this case. Now let $m
\geq 1$ be given, and assume that the result is true of all $(\theta, m)$ --
end-indiscernible sequences.

Let a formula $\phi$ and a sequence $J$ of $(\phi, m+1)$ --
end-indiscernible sequences over $\emptyset$. Let $c$ be the last element in 
$J$. Define $\psi $ so that 
\[
M\models \psi [a_0,\ldots ,a_{m-1};b]\mbox{ if and only if } M\models \phi
[a_0\cdots a_{m-1}c;b] 
\]
for all $a_0,\ldots ,a_{m-1}\in J$, $b\in M$. Note then that $|S_\psi
(B,M)\leq |S_\phi (B,M)|$ for all $B \subseteq M$, so we can use the same $
F^{*}$ for $\psi $ as for $\phi $. (This result can be improved by using a
sharper bound on the number of $\psi $ -- types.) By the definition of $g_m$
and Lemma \ref{lem7}, there must be a subset $J^{\prime \prime }$ of $J$
with cardinality at least $g_{m-1}(k-2)$ which is $\psi $ --
end-indiscernible over $A$. By the inductive hypothesis, there is a
subsequence of $J^{\prime \prime }$ with cardinality at least $k - 1$ which
is $(\psi,m) $ -- indiscernible over $A$. Form $J^{\prime }$ by adding $c$
to the end of this sequence. It follows from the $(\phi,m+1) $ --
end-indiscernibility of $J$ and the $(\psi,m) $ -- indiscernibility of $
J^{\prime \prime }$ that $J^{\prime }$ is $(\phi,m+1) $ -- indiscernible
over $A$. {\qed}

\begin{theorem}
\label{thm9} For any $A\subseteq M$ and any sequence $I$ from $M$ with $
|I|\geq g_m(k-1)$, there is a subsequence $J$ of $I$ with cardinality at
least $k$ which is $(\phi,m) $ -- indiscernible over $A$.
\end{theorem}

{\proof} By Lemmas \ref{lem7} and \ref{lem8}. {\qed}

Our goal now is to apply this to theories with different properties to see
how these properties affect the size of a sequence one must look in to be
assured of finding an indiscernible sequence. First we will do a basic
comparison between the cases when we do and do not have a polynomial bound
on the number of types over a set. In each of these cases, we will give the
bound to find a sequence indiscernible over $\emptyset $. We will use the
notation $\log ^{(i)}$ for 
\[
\begin{array}{c}
\underbrace{\log _2\circ \log _2\circ \cdots \circ \log _2} \\ 
i\mbox{ times}
\end{array}
\]

\begin{corollary}
\label{cor9}

\begin{enumerate}
\item  If $F(i)=2^{i^m}$ (which is the worst possible case), then $\log
^{(m)}g_m(k-1)\leq 4k$.

\item  If $F(i)=i^p$, then $\log ^{(m)}g_m(k-1)\leq 2m k+\log _2 k+\log _2p$.
\end{enumerate}
\end{corollary}

We now combine part (2) above with the results from the previous section to
see what happens in the specific cases of structures without the $n$ --
order property and structures without the $n$ -- independence property. We
define by induction on $i$ the function 
\[
\beth(i, x) = \left\{ 
\begin{array}{ll}
x & \mbox{if $i=0$} \\ 
2^{\beth(i - 1, x)} & \mbox{if $i>0$}
\end{array}
\right. 
\]
Recall that for the formula $\phi(x; y)$ we have defined the parameters $r=
l(x)$, $s = l(y)$, and $t = \max\{ r, s\}$.

\addtocounter{definition}{-1}

\begin{corollary}
\begin{enumerate}
\item[3.]  \setcounter{enumi}{3} If $(M,\phi )$ fails to have the $n$ --
independence property and $I=\{a_i\;:\;i<\beth (m,2k+\log _2k+\log _2n+\log
_2m)\}\subseteq M$, then there is a $J\subseteq I$ so that $|J|\geq k$ and $
J $ is a $(\phi,m) $ -- indiscernible sequence over $\emptyset $.

\item  If $(M,\phi )$ fails to have the $n$ -- order property and $
I=\{a_i\;:\;i<\beth (m,2k+\log _2k+(3ns)^{t+1})\}\subseteq M$, then there is
a $J\subseteq I$ so that $|J|\geq k$ and $J$ is a $(\phi,m) $ --
indiscernible sequence over $\emptyset $.
\end{enumerate}
\end{corollary}

Finally, note that with the additional assumption of failure of the $d$ --
cover property, if $d$ is smaller than $n$, then from the assumptions in (3)
and (4) above, we could infer a failure of the $d$ -- independence property
or the $d$ -- order property improving the bounds even further.

\section{Applications to Graph Theory}

\label{graphsec}

In this section we look to graph theory to illustrate some applications. The
reader should be warned that the word ``independent'' has a graph -
theoretic meaning, so care must be taken when reading ``independent set''
versus ``independence property''.

\subsection{The independence property in random graphs}

A first question is ``How much independence can one expect a random graph to
have?'' We will approach the answer to this question along the lines of
Albert \& Frieze \cite{af}. There an analogy is made to the Coupon Collector
Problem, and we will continue this here.

The Coupon Collector Problem (see Feller \cite{fe}) is essentially that if $
n $ distinct balls are independently and randomly distributed among $m$
labeled boxes (so each distribution has the same probability $m^{-n}$ of
occurring), then what is the probability that no box is empty? Letting $q(n,
m)$ be this probability, it is easy to compute that 
\[
q(n, m) = \sum_{i = 0}^m (-1)^i {\left( 
\begin{array}{c}
{m} \\ 
{i}
\end{array}
\right) } \left( 1 - \frac{i}{m} \right)^n = \frac{m! \, S_{m,n}}{m^n} 
\]
where $S_{n,m}$ is the Stirling number of the second kind.

It is well - known that, for $\lambda = m e^{-n/m}$, $q(n, m) - e^{-
\lambda} $ tends to 0 as $n$ and $m$ get large with $\lambda$ bounded.

The way that this will be applied in our context is as follows. We will say
that a certain set $\{v_1,\ldots ,v_k\}$ of vertices \emph{witnesses the $k$
- independence property in $G$} if $(G,R)$ has the $k$ -- independence
property with $a_i=v_i$ (see Definition \ref{def1d}). Notice that any $k$
vertices $\{v_1,\ldots ,v_k\}$ determine $2^k$ ``boxes'' defined by all
possible Boolean combinations of formulas $\{R(x,v_1),\ldots ,R(x,v_k)\}$ (a
vertex being ``in a box'' meaning it witnesses the corresponding formula in $
G$). The remaining $n-k$ vertices are then equally likely to fill each of
the $2^k$ boxes, so the probability that these $k$ vertices witness the $k$
-- independence property in $G$ is just $q(n-k,2^k)$. So for $\lambda
=\lambda (n,k)=2^k\exp (-(n-k)/2^k)$ bounded (as $n,k\rightarrow \infty $)
we will have the probability that $k$ particular vertices witness the $k$ -
independence property in a graph on $n$ vertices tends to $e^{-\lambda }$,
and the probability that a graph on $n$ vertices has the $k$ -- independence
property is at most 
\[
{\left( 
\begin{array}{c}
{r} \\ 
{k}
\end{array}
\right) }e^{-\lambda }\leq n^ke^{-\lambda }\mbox{ as }n,k\rightarrow \infty 
\]

If $n = k + k 2^k$, then $q(n-k, 2^k) \rightarrow 1$, so the particular
vertices $\{1, \ldots, k \}$ witness $k$ -- independence in a graph on $k +
k 2^k$ vertices almost surely. On the other hand,

\begin{theorem}
\label{thmg1} A random graph on $n=k+2^k(\log k)$ vertices has the failure
of the $k$ -- independence property almost surely.
\end{theorem}

{\proof} For $n=k+2^k\log k$, the $\lambda $ from above is $\frac{2^k}k$,
and $\log (n^ke^{-\lambda })=k\log (k+2^k\log k)-2^k/k$ which clearly goes
to $-\infty $ as $k\rightarrow \infty $, so the probability that a graph on $
n$ vertices has the $k$ -- independence property goes to 0. {\qed}

\subsection{Ramsey's theorem for finite hypergraphs}

We can improve (for the case of hypergraphs without $n$ -- independence) the
best known upper bounds for the Ramsey number $R_r(a, b)$. First we should
say what this means.

\begin{definition}
\begin{enumerate}
\item  An \emph{$r$ -- graph} is a set of \emph{vertices} $V$ along with a
set of $r$ -- element subsets of $V$ called \emph{edges}. The edge set will
be identified in the language by the $r$ -- ary predicate $R$.

\item  A \emph{complete $r$ -- graph} is one in which all $r$ -- element
subsets of the vertices are edges. An \emph{empty $r$ -- graph} is one in
which none of the $r$ -- element subsets of the vertices are edges.

\item  $R_r(a,b)$ denotes the smallest positive integer $N$ so that in any $
r $ -- hypergraph on $N$ vertices there will be an induced subgraph which is
either a complete $r$ -- graph on $a$ vertices or an empty $r$ -- graph on $b
$ vertices.

\item  We say that an $r$ -- graph $G$ \emph{has the $n$ -- independence
property} if $(G,R(x))$ does (where $l(x)=r$).
\end{enumerate}
\end{definition}

Note that the first suggested improvement of Lemma \ref{lem7} applies in
this situation -- namely, the edge relation is symmetric. We can immediately
make the following computations.

\begin{lemma}
\label{lemg2}

\begin{enumerate}
\item  In an $r$ -- graph $G$, $F$ is given by $F(i)=2^q$ where $q={\left( 
\begin{array}{c}
{i} \\ 
{r-1}
\end{array}
\right) }$. Consequently, $F^{*}(k)\leq 2^{k^r}$ in this case.

\item  In an $r$ -- graph $G$ which does not have the $n$ -- independence
property, $F$ is defined by 
\[
F(i):=\left\{ 
\begin{array}{ll}
1 & \mbox{for }i<r \\ 
i^{(r-1)(n-1)} & \mbox{otherwise}
\end{array}
\right. 
\]
Consequently $F^{*}(k)\leq k^{(r-1)(n-1)k}$ in this case.
\end{enumerate}
\end{lemma}

For a fixed natural number $p$, define the functions $E_p ^{(j)}$ by

\begin{itemize}
\item  $E^{(1)}=E=(\alpha \mapsto (\alpha +1)^{p(\alpha +1)})$, and

\item  $E^{(i+1)}=E\circ E^{(i)}$ for $i\geq 1$.
\end{itemize}

\begin{theorem}
\label{thmg2} Let $n \geq 2$ and $k \geq 3$ be given, and let $p = (r-1)(n-1)
$. If an $r$ -- graph $G$ on at least $E_p ^{(r-1)}(k-1)$ vertices does not
have the $n$ -- independence property, then $G$ has an induced subgraph on $k
$ vertices which is either complete or empty.
\end{theorem}

{\proof} (By induction on $r$)

For $r=2$, the graph has at least $E_{n-1}^{(1)}(k-1) = k^{(n-1)k}\geq 2^{2k}
$ vertices, and it is well-known (see e.g. \cite{grs}) that $
2^{2k}\rightarrow (k)_2^2$.

Let $r \geq 3$ be given, and let $G=(V,R)$ be an $r$ -- graph as described
and set $C = E_p ^{(r-2)}(k)$, where $p = (r-1)(n-1)$. Using $F(i)=i^p$ for $
i\geq 2$, ($F(0)=F(1)=1$) and computing $F^{*}$ in Lemma \ref{lemg2}, we
first see from Lemma \ref{lem7} that any $r$ -- graph on at least $
(C+1)^{p(C+1)}$ vertices will have an $(R, 1)$ -- end-indiscernible sequence 
$J$ over $\emptyset $ of cardinality $C$. Let $v$ be the last vertex in $J$
and define the relation $R^{\prime}$ on the $(r-1)$ -- sets from (the range
of) $J$ by

\begin{center}
$R^{\prime }(X)$ if and only if $R(X \cup \{v\})$.
\end{center}

Now $(J,R)$ is an $(r-1)$ -- graph of cardinality $C$, so by the inductive
hypothesis there is an $R^{\prime }$ -- indiscernible subsequence $J_0$ of $J
$ with cardinality $k$. Clearly $I=\{A \cup \{v\} : A \in J_0\}$ is an $R$
-- indiscernible sequence over $\emptyset $ of cardinality $k$. {\qed}

\noindent\textsc{Remark:} Another way to say this is that in the class of $r$
-- graphs without the independence property $R_r(k,k) \leq E_p ^{(r-1)}(k-1)$
.

\subsection*{Comparing upper bounds for $r=3$}

Note that for $r=3$ in Theorem \ref{thmg2}, we have $p=2(n-1)$, and so we
get $E_p^(2)(k-1) = (2^{2k}+1)^{p(2^{2k}+1)}$ which is roughly $
2^{nk(2^{2k}+2)}$. The upper bound for $R_3(k,k)$ in \cite{ehmr} is roughly $
2^{2^{4k}}$. So $\log _2\log _2($their bound$)=4k$ and 
\[
\log _2\log _2(\mbox{our bound})=\log _2(nk(2^{2k+2}))=\log _2p+\log
_2k+(2k+2) 
\]
which is smaller than $4k$ as long as $2k-2-\log _2k>\log _2n$. This is true
as long as $n<2^{2k-2}/k$.

For example, for $k = 10$ our bound is about $2^{c(n-1)}$ where $c$ is
roughly $4 \times 10^7$ and theirs is about $2^{2^{40}}$. Since $2^{40}$ is
roughly $10^{12}$, this is a significant improvement in the exponent for 3
-- graphs without the $n$ -- independence property.

\subsection*{Comparing upper bounds in general}

Let $a_r$ be the upper bound for $R_r(k,k)$ given in \cite{ehmr} and $b_r$
be the upper bound as computed for the class of $r$ -- graphs without the $n$
-- independence property in Theorem \ref{thmg2} (both as a function of $k$,
the size of the desired indiscernible set). Since we have $b_{r+1} \leq
b_r^{(p)(b_r)}$, we get the relationship 
\begin{eqnarray*}
\log^{(r)} b_{r+1} & \leq & \log^{(r-1)} [p \, b_r (\log b_r)] \\
& = & \log^{(r-2)} (\log (r-1) + \log (n-1) + \log b_r + \log \log b_r)
\end{eqnarray*}
for $r \geq 3$, $\log \log b_3 = 2k + \log_2 k + \log_2 n + \log_2 r$, and $
\log b_2 = 2k$. It follows that $\log^{(r)} b_{r+1}$ is less than (roughly) $
2k + \log_2 k + \log_2 n$ for every $r$.

In \cite{ehmr}, the bounds $a_r$ satisfy $\log a_2 = 2k$, $\log \log a_3 =
4k $, and for $r \geq 3$, 
\begin{eqnarray*}
\log^{(r)} a_{r+1} & = & \log^{(r-1)} (a_r^r) = \\
\log^{(r-2)}(r \log a_r) & = & \log^{(r-3)} [\log r + \log \log a_r].
\end{eqnarray*}
We can then show that $\log^{(r-1)} a_r < 4k + 2$ for all $r$.

Clearly for each $r \geq 3$, $\frac{b_r}{a_r} \rightarrow 0$ as $m$ gets
large.

\noindent\textsc{Final Remark:} On a final note, the above comparison is
only given for $r$ -- graphs with $r \geq 3$ because the technique enlisted
does not give an improvement in the case of graphs. This has not been
pursued in this paper because it seems to be of no interest in the general
study. However, the techniques may be of interest to the specialist.

\section{Toward a classification theory}

\subsection{Introduction}

One of the most powerful concepts of model theory (discovered by Shelah) is
the notion of forking, which from a certain point of view can be considered
as an instrument to discover the structure of combinatorial geometry in
certain definable subsets of models.

We were unable to capture the notion of forking (or a forking - like
concept) for finite structures. What we \emph{can} do is to present an
alternative, more global property called stable -- amalgamation (which is
the main innovation in \cite{sh3} in dealing with non-elementary classes).
We do this by imitating \cite{sh4}. We have reasonable substitutes for $
\kappa(T)$ and $Av(I,A,M)$. The most important property we manage to prove
is the symmetry property for stable amalgamation. This is the corresponding
property to the exchange principle in combinatorial geometry.

The ultimate goal of the project started here is to have a decomposition
theorem not unlike the theorem for finite abelian groups. We hope to
identify some properties $P_1, \ldots, P_n$ of a class of finite models $K$
in such a way that the following conjecture will hold:

\begin{conjecture}
\label{dream} If $\langle K,\prec_K \rangle$ satisfies $P_1,\ldots,P_n$ then
for every $M\in K$ large enough there exists a finite tree $T\subseteq
\;{}^{\omega <}\omega$ and $\{M_{\eta} \prec_K M\; : \; \eta\in T \}$ such
that

\begin{enumerate}
\item  $\{M_{\eta} \; : \; \eta \in T \}$ is a ``stable'' tree \footnote{
Defined using the notion of stable amalgamation introduced below.}

\item  For every $\eta\in T$ we have that $\| M_{\eta} \| \leq n(K)$ 
\footnote{ 
We acknowledge that cardinality of the universe may not be an appropriate
measure of ``smallness'' for a substructure in this context. The reader
should also consider a restriction on the cardinality of a set of generators
for $M_{\eta}$ as another possibility.}

\item  $M$ is uniquely determined by $\bigcup_{\eta\in T}|M_{\eta}|$.
\end{enumerate}
\end{conjecture}

Ideally $P_1,\ldots,P_n$ is a minimal list of properties sufficient to
derive the above decomposition. We hope to be able to eventually emulate
Theorem XI.2.17 in \cite{sh}.

Considering our present state of knowledge it seems that our conjecture is
closer to a fantasy than to a mathematical statement. However we seem to
have a start. Much of this section together with some of the earlier results
can be viewed as a search for candidates for the above mentioned list of
properties $P_1,\ldots,P_n$. We hope that this section might form an
infrastructure for the classification project.

\subsection{Abstract properties}

We now begin to look at some of the abstract properties of a class $K$ of
finite $L$ -- structures with an appropriate partial ordering denoted by $
\prec_K$. These properties come from Shelah's list of axioms in \S 1 of \cite
{sh3}.

\begin{definition}
Let $L$ be a given similarity type, let $\Delta $ be a set of $L$ --
formulas, and let $n<\omega $, by $\Delta _n^{*}$ we denote the minimal set
of $L$ -- formulas containing the following set and all its subformulas: 
\[
\{\exists x[\bigwedge_{i\in w}\phi (x;y_i)\wedge \bigwedge_{i\in k \setminus
w} \neg \phi (x;y_i)]:\phi (x;y)\in \Delta , k \leq n, w\subseteq
k,l(y_i)=l(y)\}. 
\]
\end{definition}

We will now look into natural values of $k$ from the previous section.

\begin{theorem}
Suppose the formula $\phi $ fails to have the weak $n$ -- order property
(and hence fails to have the $n$ -- independence property in $K$), and let $
\Delta \supseteq \{ \phi \}^*_n$ be given.  Then for every $(\Delta,n)$ --
indiscernible sequence $I$ over $\emptyset$ and every $c \in M$, either $
|\{a \in I \; : \; M \models \phi [c;a] \}| < n$ or $|\{a \in I \; : \; M
\models \neg \phi [c;a] \}| < n$.
\end{theorem}

{\proof} We may assume the length of $I$ is at least $2n$ since otherwise
the result is trivial. We proceed by contradiction. Suppose the result is
not true. Then there is a $(\Delta,n)$ -- indiscernible sequence $I$ from $M$
of length at least $2n$ and a $c\in M$ such that both $|\{\phi
[c;a]\;:\;a\in I\}|\geq n$ and $|\{\neg \phi [c;a]\;:\;a\in I\}|\geq n$. Let 
$\{a_0,\ldots ,a_{2n-1}\}\subseteq I$ be such that 
\begin{equation}
M\models \bigwedge_{i<n}\phi [c,a_i]\;\wedge \;\bigwedge_{n \leq i < 2n}
\neg \phi [c,a_i]  \label{num1}
\end{equation}
We complete the proof by showing that $\{a_0,\ldots ,a_{n-1}\}$ exemplifies
the $n$ -- independence property. Let $w\subseteq n$ be given. Consider the
formula 
\[
\psi _w(y_0,\ldots ,y_{n-1})\stackrel{\mathrm{def}}{=}(\exists x)\left[
\bigwedge_{i\in w}\phi (x;y_i)\;\wedge \;\bigwedge_{i\in n \setminus w}\neg
\phi (x;y_i) \right] 
\]
Let $\{i_0,\ldots ,i_{k-1}\}$ be an increasing enumeration of $w$. By (\ref
{num1}) the following holds 
\[
M\models \psi _w[a_{i_0},\ldots ,a_{i_{k-1}},a_n,\ldots ,a_{2n-1-k}] 
\]
Since $\psi _w\in \{\Delta \}^*_n$, by the indiscernibility of $I$ we have
that also $M\models \psi _w[a_0,\ldots ,a_{n-1}]$. So for every $w\subseteq n
$, we may choose $b_w\in M$ so that 
\[
M\models \bigwedge_{i\in w}\phi [b_w,a_i]\;\wedge \;\bigwedge_{i\in n-w}\neg
\phi [b_w,a_i]. 
\]
We are done since $\{a_0,\ldots ,a_{n-1}\}$ and $\{b_w:w\subseteq n\}$
witness the fact that $(M,\phi )$ has the $n$ -- independence property. {
\qed}

The following definition is inspired by $\kappa (T)$ in Chapter III of \cite
{sh}.

\begin{definition}
Let $n<\omega $ be given, and let $\Delta $ be a finite set of formulas. $
\kappa _{\Delta ,n}(K)$ is the least positive integer so that for every $
M\in K$, every sequence $I=\langle a_i\;:\;i<\beta <\omega \rangle \in M$
which is $\Delta _n^{*}$ -- indiscernible over $\emptyset$ has either $
M\models \phi [c;a_i]$ or $M\models \neg \phi [c;a_i]$ for less than $\kappa
_{\Delta ,n}(K)$ elements of $I$ for each $\phi \in \Delta $ and $c\in M$.
Recalling that $\Delta$ is to be closed under negation, we will write $
\kappa_{\phi, n}$ instead of $\kappa_{ \{ \phi, \neg \phi \},n}$.
\end{definition}

So the previous theorem states that if the formula $\phi$ fails to have the $
n$ -- independence property in $M$, then $\kappa_{\phi,n}(M) \leq n$. When $
\phi$ is understood to  not have the $n$ -- independence property $
\kappa_{\phi}$ will stand for $\kappa_{\phi,n}$. In this case, the following
definition makes sense.

\begin{definition}
Let $n < \omega$ be given, $\Delta $ be a finite set of formulas and $n$ be
as above. Suppose $I$ is a sequence of $(\Delta_n^{*},n)$-indiscernibles
over $\emptyset $. Define 
\begin{eqnarray*}
Av_\Delta (I,A,M)& = & \{\phi (x;a):a\in A,\phi (x;y)\in \Delta \\
& & \mbox{ and } |\{c\in I:M\models \phi [c;a]\}|\geq \kappa _{\Delta,n}
(K)\}.
\end{eqnarray*}
If $M$ is understood, it will often be omitted.
\end{definition}

\begin{theorem}
Let $\psi(y; x) = \phi(x; y)$. If $\phi $ has neither the $n$ --
independence property nor the $d$ -- cover property in $M$, $\Delta
\supseteq \{\phi \}_n^{*}$ is finite, and $I$ is a set of $(\Delta,n) $ --
indiscernibles over $\emptyset $ of length greater than $\max \{d \cdot
\kappa _\psi (M),2n\}$, then $Av_\phi (I,A,M)$ is a complete $\phi $ -- type
over $A$.
\end{theorem}

{\proof} That $Av_\psi (I,A,M)$ is complete follows from the previous
theorem. To see that it is consistent, we need only establish that every $d$
formulas from it are consistent (by failure of the $d$ -- cover property),
and this follows >from the size of $I$ and the pigeonhole principle. {
\qed}

So we will use the following term to denote when we are in a model in which
the notion of average type is well defined.

\begin{definition}
Let $\psi(y; x) = \phi(x; y)$, and $\Delta = \{ \phi, \psi, \neg\phi,
\neg\psi \}$. We will say that $M$ is \underline{$(\phi ,n,d)$ -- good} if $
(M,\Delta )$ has neither the $n$ -- independence property nor the $d$ --
cover property. In this case, we will define $\lambda _\phi (M)= \max \{d
\cdot \kappa _{\Delta,n}(M),2n \}$. We will sometimes refer to this same
situation by saying $(M, \phi)$ is $(n,d)$ -- good.

If $K$ is a class of $(\phi, n, d)$ -- good structures which all include a
common set $A$, then we will say that $K$ is $(\phi, n, d)$ -- good, and
define $\lambda (K)= \kappa _{\Delta,n}(K) \cdot |A|^s$, where $s = \max
\{l(y) : \phi(x; y) \in \Delta \}$.
\end{definition}

\begin{example}
\label{ex.fields}Let $T$ be an $\aleph _1$ -- categorical theory in a
relational language (no function symbols), and let $M\models T$ be an
uncountable model (e.g., an uncountable algebraically closed field of
positive characteristic). Let $K:=\{N\subseteq M:\Vert N\Vert <\aleph _0\}$,
and let $\phi \in L(T)$. By $\aleph _1$ -- categoricity there exist integers 
$n$ and $d$ such that $K$ is $(\phi, n,d)$ -- good (see Corollary \ref
{finitecor}).
\end{example}

\begin{definition}
\label{submodel} For a fixed (finite) relational language $L$ and an $L$ --
formula $\phi $ (and $\psi(y; x) = \phi(x; y)$), let $K$ be a class of
finite $(\phi ,n,d)$ -- good $L$ -- structures all of which include a common
set $A$. Fix a $k < \omega$ and and define $\prec _{K}$, as follows: $N\prec
_{K} M$

\begin{enumerate}
\item  $N \subseteq M$, and for all $a \in A$ and $b \in N$, $M \models \phi[
b; a]$ if and only if $N \models \phi[b; a]$.

\item  For every $a_0, \ldots, a_{k-1} \in A$, if $M\models \exists
x\bigwedge_{i<k}\phi (x;a_i)$, then $M\models \bigwedge_{i<k} \phi [b;a_i]$
for some $b\in N$.

\item  For every $a\in M$, there is a sequence $I \subseteq N$ which is $
(\{\psi \}^*_n, n)$ -- indiscernible over $A$ with length at least $\lambda
(K)$ so that $tp_\phi (a,A,M)=Av_\phi(I,A,M)$.
\end{enumerate}

We define the same relation for a set $\Delta $ of formulas simply by
requiring that the above holds for each $\phi \in \Delta$ in the case that $K
$ is a class of finite $(\phi ,n,d)$ -- good structures.
\end{definition}

\noindent\textsc{Remarks on Definition \ref{submodel}:}

\begin{itemize}
\item  Condition (1) ensures that the fact for elementary classes that forms
the basis of the Tarski - Vaught (namely, $N \subseteq M \Rightarrow N
\prec_{qf} M$ holds for $\phi$ -- formulas even if $\phi$ is not quantifier
free.

\item  Condition (2) is like $k$ -- saturation relative to $\phi $ --
formulas with parameters from $A$ (i.e. every $\phi$ -- type with at most $k$
parameters from $A$ which is realized in $M$ is also realized in $N$). It
can be thought of as a generalization of the Tarski - Vaught test
relativized to formulas from $\{ \phi \}^*_k$.

\item  Condition (3) is a property like the one that guarantees in the first
order case that types over a model are stationary. Here we are requiring a
strong closure condition on $N$ --- namely, if $a \in M \setminus N$, then
there is a strong reason why $a$ does not belong to $N$ : there is a long
sequence of indiscernibles in $N$ averaging the same $\phi $ -- type over $A$
.

\item  It should be emphasized that $K$ (and hence $\prec_K$) has parameters 
$A$, $k$, and $\phi$ which are suppressed for notational convenience.
\end{itemize}

\subsection{Properties of $\prec_K$}

We prove the following facts about the relation $\prec_{K}$. The Roman
numerals in parentheses indicate the corresponding Axioms in \cite{sh3}. $
\Delta$ is closed under negation and fixed throughout, and $K$ is a class of 
$(\Delta, n, d)$ -- good structures which all include a common set $A$.

\begin{lemma}
\begin{enumerate}
\item  (I) If $N\prec _K M$, then $N\subseteq M$.

\item  (II) $M_0 \prec_K M_1 \prec_K M_2$ implies $M_0 \prec_K M_2$. Also $
M\prec _KM$ for all $M \in K$

\item  (V) If $N_0\subseteq N_1\prec _KM$ and $N_0\prec _KM$, then $N_0\prec
_KN_1$.
\end{enumerate}
\end{lemma}

{\proof} The first part of (I) is trivial, and the second part of (II) only
requires that one chooses a constant sequence for $I$ in Condition (3). Note
that for all three statements, checking Condition (1) is routine.

For the first part of (II), assume the hypothesis is true and first look at
Condition (2). Let $\phi \in \Delta$ and $a_i \in A$ for $i<k$ be given, and
assume that $M_2 \models \exists x\bigwedge_{i<k}\phi (x;a_i)$ Since $M_1
\prec_K M_2$, we may choose $b \in M_1$ so that $M_2 \models \bigwedge_{i<k}
\phi [b;a_i]$. Of course, $b, a_0, \ldots, a_{k-1}$ are all from $M_1$, so
we can conclude that $M_1 \models \bigwedge_{i<k}\phi [b;a_i]$, or less
specifically that $M_1 \models \exists x\bigwedge_{i<k}\phi (x;a_i)$. Since $
M_0 \prec_K M_1$, this in turn allows us to choose $b^{\prime} \in M_0$ so
that $M_1 \models \bigwedge_{i<k}\phi [b^{\prime};a_i]$, which means
necessarily that $M_2 \models \bigwedge_{i<k}\phi [b^{\prime};a_i]$.

For condition (3), let $a \in M_2$ be given. Since $M_1 \prec_K M_2$, we may
choose $I$ from $M_1$ of length $\lambda(K)$ which is $(\{\psi \}^*_n, n)$
-- indiscernible over $A$ so that $tp_\phi (a,A,M_2) = Av_\phi(I,A,M_2)$.
Because the length of $I$ exceeds $|A|^{l(y)} \cdot \kappa_{\phi}(M)$, we
may choose one element $b_{i_0}$ in $I$ so that $tp_{\phi}(b_{i_0}, A, M_2)
= Av_{\phi}(I, A, M_2)$. (This can be accomplished by throwing out $<
\kappa_{\phi}(M)$ elements of $I$ for each instance $\phi(x;b)$ with $b \in A
$ so that the elements of $I$ that are left all realize the same instances
of $\phi$ over $A$.) Since $M_0 \prec_K M_1$, we may choose a sequence $J$
in $M_0$ which is $(\{\psi \}^*_n, n)$ -- indiscernible over $\emptyset$ so
that $Av_{\phi}(J, A, M_1) = tp_{\phi}(b_0, A, M_1)$. But then we have $
Av_{\phi}(J, A, M_1) = Av_{\phi}(I, A, M_1)$, and consequently $Av_{\phi}(J,
A, M_2) = Av_{\phi}(I, A, M_2) = tp_{\phi}(a,A,M_2)$, as desired.

For (V), consider first Condition (2). Assuming the hypotheses in (V) are
true, we let $\phi \in \Delta$ and $a_i \in A$ for $i<k$ be given, and
assume that $N_1 \models \exists x\bigwedge_{i<k}\phi (x;a_i)$. It then
follows from $N_1 \prec_K M$ (Condition (2)) that $M \models \exists
x\bigwedge_{i<k}\phi (x;a_i)$, and since $N_0 \prec_K M$, we may choose $b
\in N_0$ so that $M\models \bigwedge_{i<k} \phi [b;a_i]$. Of course, $b,
a_0, \ldots, a_{k-1}$ are all from $N_1$, so >from Condition (1) of $N_1
\prec_K M$, we can conclude that $N_1 \models \bigwedge_{i<k}\phi [b;a_i]$.

For condition (3), let $a \in N_1$ be given. Since $a \in M$ and $N_0
\prec_K M$, we may choose $I$ from $N_0$ of length $\lambda(K)$ which is $
(\{\psi \}^*_n, n)$ -- indiscernible over $\emptyset$ so that $tp_\phi
(a,A,M) = Av_\phi(I,A,M)$. Since $N_1 \subseteq M$, and $a$, $I$, and $A$
are all included in $N_1$, it follows that $tp_\phi
(a,A,N_1)=Av_\phi(I,A,N_1)$. {\qed}

\begin{definition}
Here again $\psi (y;x)=\phi (x;y)$. Given $(\phi ,n,d)$ -- good structures $M
$, $M_0$, $M_1$, and $M_2$ with $M_l\prec _KM$, $M_0\prec _KM_1$, and $
M_0\prec _KM_2$, we say that $(M_0,M_1,M_2)$ is in $\phi $ -- \underline{
stable amalgamation} inside $M$ if for every $c\in M_2$ with $l(c)=l(x)$
there is a $(\{\psi \}_n^{*},n)$ -- indiscernible sequence $I\subseteq M_0$
over $\emptyset $, of length at least $\lambda (K)$ such that $Av_\phi
(I,M_1,M)=tp_\phi (c,M_1,M)$.
\end{definition}

To prove symmetry of stable amalgamation (with the assumption of non-order),
we must first establish the following lemma (corresponding to I.3.1 in \cite
{sh4}).

\begin{lemma}
Let $\psi (y;x)=\phi (x;y)$ and $\Delta =\{\phi ,\psi ,\neg \phi ,\neg \psi
\}$, and let $\lambda =\max \{\lambda _\Delta (M),\kappa _\phi (M)+\kappa
_\psi (M)+\kappa _\phi (M)\cdot \kappa _\psi (M)\}$. Assume $M$ is $(\Delta
,n,d)$ -- good, $I_0=\langle a_k^0:k<m_0\rangle $ is a $\{\psi \}_n^{*}$ --
indiscernible sequence (over $\emptyset $) in $M$ of length greater than $
\lambda $, and $I_1=\langle a_k^1:k<m_1\rangle $ is a $\{\phi \}_n^{*}$ --
indiscernible sequence (over $\emptyset $) in $M$ of length greater than $
\lambda $. The following are equivalent:

\begin{enumerate}
\item[(i)]  There exists $i_k<m_0$ for $k<m_0-\kappa _\phi (M)$ such that
for each $k$, 
\[
\phi (a_{i_k}^0,y)\in Av_\psi (I_1,|M|,M).
\]

\item[(ii)]  There exists $j_l<m_1$ for $l<m_1-\kappa _\psi (M)$ such that
for each $l$, 
\[
\phi (x,a_{j_l}^1)\in Av_\phi (I_0,|M|,M).
\]
\end{enumerate}
\end{lemma}

{\proof} Assume (i) holds. Choose $i_k<m_0$ for $k<m_0-\kappa _{\phi }(M)$
and $j_{k,l}<m_1$ for $l<m_1-\kappa _{\psi }(M)$ witnessing (i). Since $m_0
> \kappa _{\psi }(M)+\kappa _{\phi }(M)\kappa _{\psi }(M)$, we can find $
\kappa _{\psi }(M)$ of the $j_{k,l}$ each of which occurs for at least $
\kappa _\phi (M)$ different $i_k$. Thus for each of these, $
\phi(x;a_{j_{k,l}}^1)\in Av_\phi (I_0,|M|,M)$.

Now assume (ii) does not hold. That is, there are $j_l<m_1$ for each $
l<m_1-\kappa _{\psi }(M)$ such that $\neg \phi (x;a_{j_l}^1)\in Av_\phi
(I_0,|M|,M)$. Clearly one of these $j_l$ must correspond to one of the $
j_{k,l}$ from before that occurs at least $\kappa _{\psi }(M)$ times. But as
we noted above $\phi (x;a_{j_{k,l}}^1)\in Av_\phi (I_0,|M|,M)$, a
contradiction.

Note that (ii) implies (i) by the symmetric argument. {\qed}

\begin{theorem}[Symmetry]
Let $\psi (y;x)=\phi (x;y)$ and $\Delta =\{\phi ,\psi ,\neg \phi ,\neg \psi
\}$. Let $K$ be a class of $(\Delta ,n,d)$ -- good structures which all
include a common set $A$. Suppose $M_0$, $M_1$, $M_2\prec _KM$, $M_0\prec
_KM_1$, and $M_0\prec _KM_2$. Then $(M_0,M_1,M_2)$ is in $\Delta $ -- stable
amalgamation inside $M$ if and only if $(M_0,M_2,M_1)$ is in $\Delta $ --
stable amalgamation inside $M$.
\end{theorem}

{\proof} We show that $(M_0, M_1, M_2)$ in $\phi $ -- stable amalgamation
implies that $(M_0,M_2,M_1)$ is in $\psi $ -- stable amalgamation. The
result follows from this. Assume that $(M_0,M_1,M_2)$ is in $\phi $ --
stable amalgamation in $M$. Let $c\in M_1$ with $l(c)=l(x)$ be given. (We
need to find a $(\{\phi\}^*_n, n) $ -- indiscernible sequence $I\subseteq M$
, with $Av_\psi (I,M_2,M)=tp_\psi(c,M_2,M)$.) By the definition of $M_0\prec
_KM_1$, we may choose a $(\{\phi\}^*_n, n)$ -- indiscernible $I\subseteq M_0$
of length at least $\lambda (K)$ so that $tp_\psi (c,M_0,M_1)=Av_\psi
(I,M_0,M_1)$ (and so $tp_\psi (c,M_0,M)=Av_\psi (I,M_0,M)$).

We claim that $Av_\psi (I,M_2,M)=tp_\psi (c,M_2,M)$. (Note that the first
type is defined since $I$ is long enough.) To see this, let $b\in M_2$ be
given such that $M\models \psi [c;b]$, and we will show that $\psi (x;b)\in
Av_\psi (I,M_2,M)$.

Since $(M_0,M_1,M_2)$ is in $\phi $ -- stable amalgamation in $M$, we may
choose a $(\{\psi\}^*_n, n)$ -- indiscernible set $J\subseteq M_0$ of length
at least $\lambda (K)$ so that $tp_\phi (b,M_1,M)=Av_\phi (J,M_1,M)$. Since $
M\models \phi [b;c]$, we have $\phi (x;c)\in Av_\phi (J,M_1,M)$, so a large
number of $b_i$ from $J$ have $M\models \phi [b_i;c]$, or rather $\psi
(y;b_i)\in tp_\psi (c,M_0,M)=Av_\psi (I,M_0,M)$ for each of these $b_i$. So $
\psi (y;b_i)\in Av_\psi(I,M_0,M)$ for each of these $b_i$.

But then by the previous Lemma, we may choose a large number of $c_j$ from $I
$ for which $\phi (x;c_j)\in Av_{\phi}(J,M_1,M) = tp_{\phi}(b, M_1, M)$.
That is, $M\models \phi[b;c_j]$ for each of these $c_j$, and so $\psi
(y;b)\in Av_\psi (I,M_2,M)$ as desired. {\qed}

\end{document}